\pdfoutput=1
\documentclass{amsart}
\usepackage{amscd}
\usepackage{amssymb}
\usepackage{graphicx}
\usepackage{hyperref}

\theoremstyle{plain}
\newtheorem{theorem}{Theorem}

\newtheorem{lemma}{Lemma}
\newtheorem{fact}{Fact}

\newtheorem{question}{Question}

\theoremstyle{definition}

\theoremstyle{remark}

\begin{document}
\title[Singular Legendrian knots]{On elementary moves of singular Legendrian knots}
\author{Sara Yamaguchi}
\email{st20190sr@gm.ibaraki-ct.ac.jp}
\address{National Institute of Technology, Ibaraki College, 312-8508, Japan}
\author{Noboru Ito}
\email{nito@gm.ibaraki-ct.ac.jp}
\address{National Institute of Technology, Ibaraki College, 312-8508, Japan}
\keywords{Legendrian front; Legendrian knot; contact structure; Legendrian singular knot; singular knot; plane curve}
\date{January 11, 2022}
\maketitle

\begin{abstract}
We have two results.  
First, we give $96$ generating sets oriented singular Reidemeister moves; it is an answer to a question by Bataineh, Khaled, Elhamdadi, and Hajij who give a generating set of oriented singular Reidemeister moves using their computation.     
Second, in the theory of plane curve and Legendrian knots introduced by V.~I.~Arnold, we select which moves survive as those of Legendrian singular knots and fronts diagrammatically and explicitly.    
\end{abstract}

\section{Introduction}\label{intro}
The standard contact structure of $\mathbb{R}^3$ with the contact $1$-form $dz+xdy$ has been well studied, e.g., \emph{Legendrian Reidemeister theorem} is known;   
on the other hand, plane curves with Arnold's \cite{Arnold1994} contact $1$-form is $\cos \theta dx+ \sin\theta dy$ also implying Legendrian knots in the solid torus \cite{ChmutovGoryunov1997, HayanoIto2015}.    Nowadays Arnold's approach has been also studied as  an analogue of Vassiliev theory of knots, e.g., Vassiliev-type skein relation has been introduced \cite{Goryunov1998}.   However, as far as we know, there is few study toward giving Legendrian Reidemeister theorem or even candidates of necessary \emph{diagrammatic generating sets}.  Therefore, in this paper, we give a sets of diagrammatic moves of Legendrian singular knots and fronts derived from those of singular knots.  

\begin{figure}
\includegraphics[width=8cm]{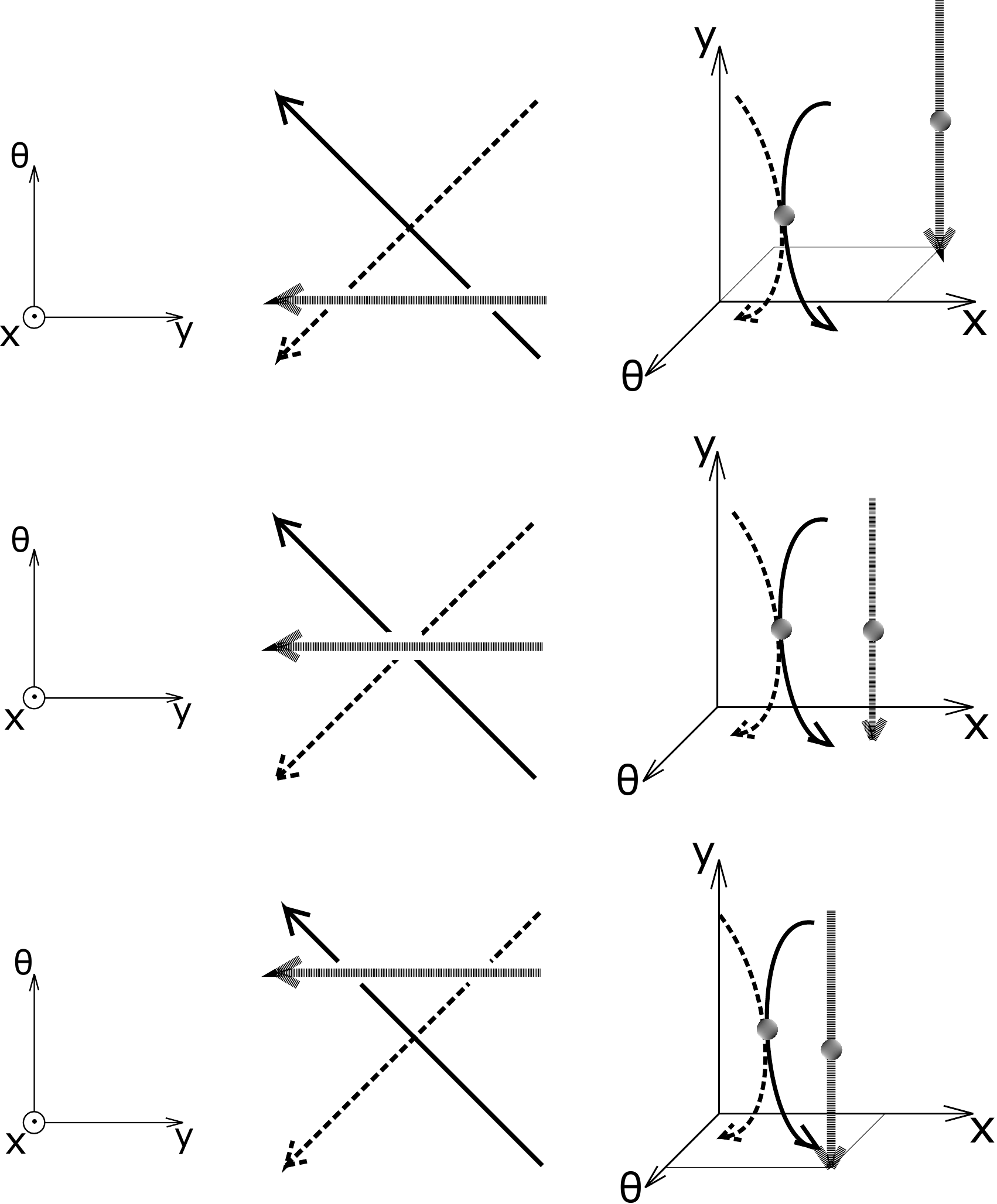}
\caption{Front presentation of $\Omega_{4b}$.  }\label{ronbunL}
\end{figure}

\begin{figure}
\includegraphics[width=8cm]{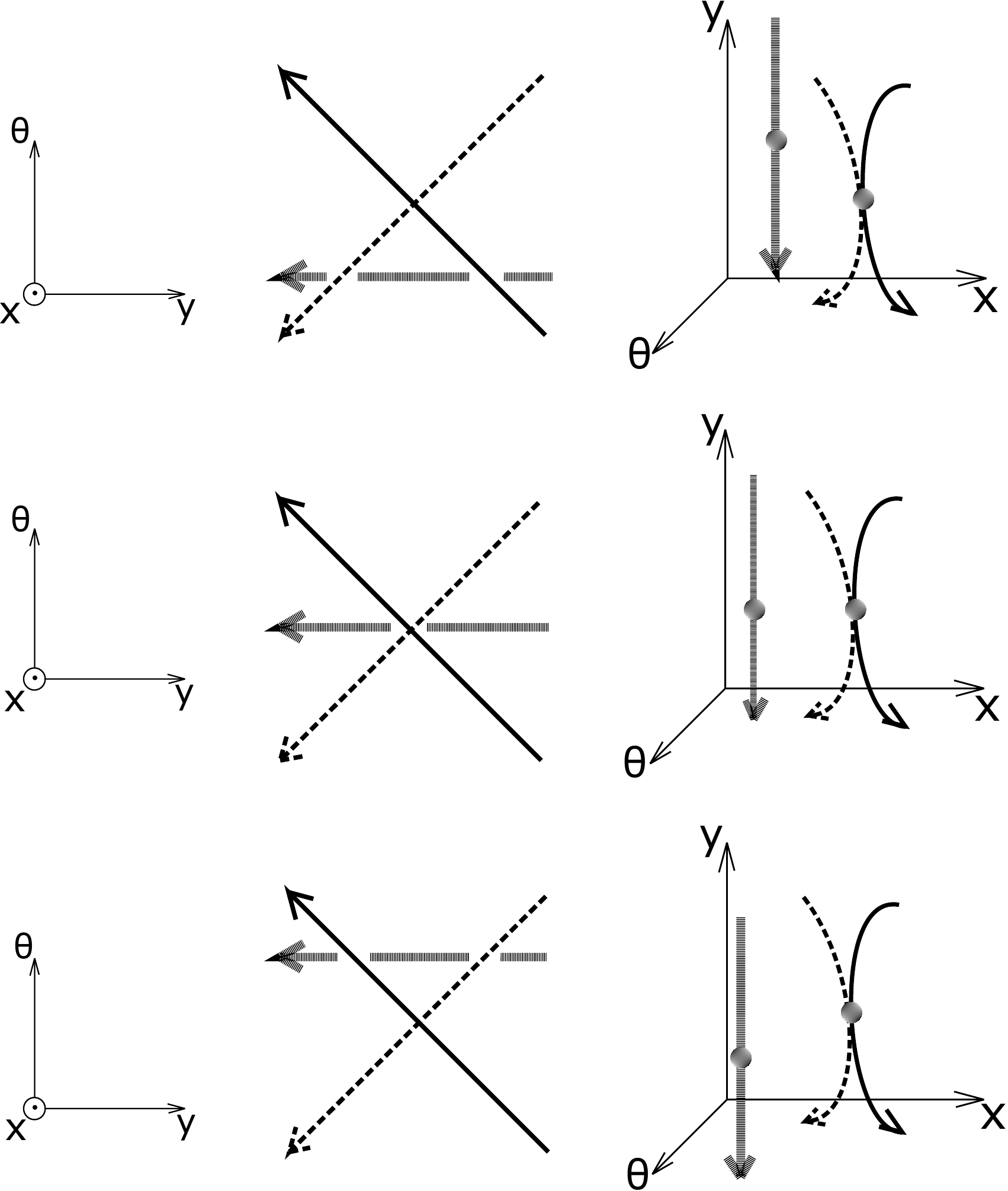}
\caption{Front presentation of $\Omega_{4f}$.  }\label{4f}
\end{figure}

\begin{figure}
\centering
\includegraphics[width=8cm]{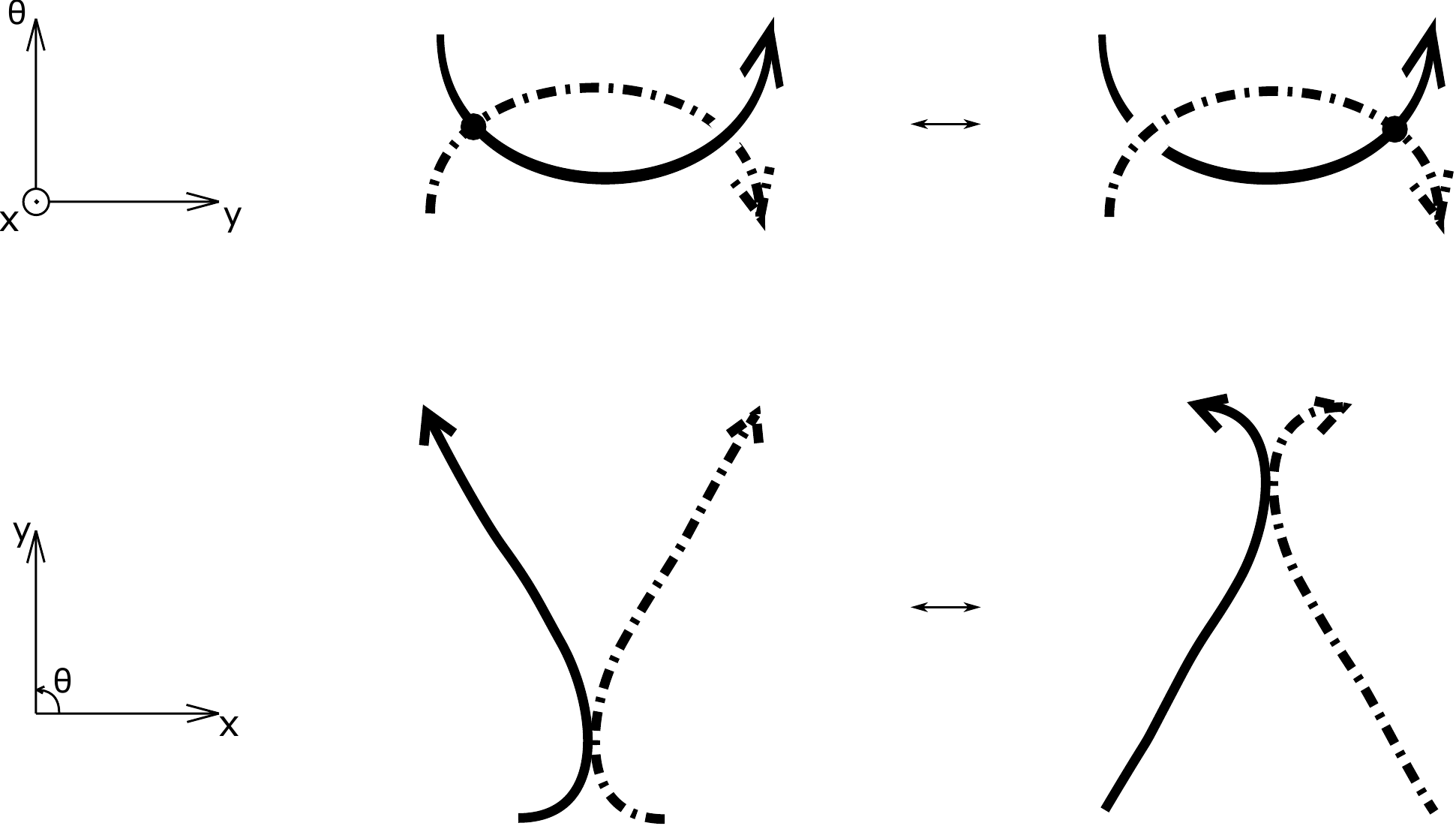}
\caption{Front presentation of $\Omega_{5d}$.  Other $\Omega_{5*}$ does not have a front presentation.}\label{ronbunN}
\end{figure}

In this paper, symbols of singular Reidemeister moves obey those of Bataineh, Khaled, Elhamdadi, and Hajij \cite{BatainehKhaledElhamdadiHajij2018}.   They give Facts~\ref{thmBKEHa} and \ref{thmBKEHb} up to nonsingular Reidemeister moves and call Rediemeister moves with respect to singular points \emph{singular Reidemeister moves}.  
\begin{fact}[\cite{BatainehKhaledElhamdadiHajij2018}]\label{thmBKEHa}
Only three oriented singular Reidemeister moves are required to generate the entire set of $\Omega_4$, and $\Omega_5$ moves.  These three moves are $\Omega_{4a}$, $\Omega_{4e}$, and $\Omega_{5a}$.   
\end{fact}
\begin{fact}[\cite{BatainehKhaledElhamdadiHajij2018}]\label{thmBKEHb}
Only two oriented singular Reidemeister moves of type $\Omega_4$ are required to generate all type $\Omega_4$ moves.  These moves are $\Omega_{4a}$ and $\Omega_{4e}$. 
\end{fact}
In \cite{BatainehKhaledElhamdadiHajij2018}, they gave open questions.  One of them is Question~\ref{qGen}.  
\begin{question}\label{qGen}
Find other generating sets of oriented singular Reidemeister moves.  
\end{question}
For Question~\ref{qGen}, we have: 
\begin{theorem}\label{thmGenSing}
Generating sets of singular Reidemeister moves are listed by $96$ cases: any triple $x, y, z$ such that $x$ in $\{ \Omega_{4a}, \Omega_{4b}, \Omega_{4c}, \Omega_{4d} \}$, $y$ in $\{ \Omega_{4e}, \Omega_{4f}, \Omega_{4g}, \Omega_{4h} \}$, and $z$ in $\{ \Omega_{5a}, \Omega_{5b}, \Omega_{5c}, \Omega_{5d}, \Omega_{5e}, \Omega_{5f} \}$.  
\end{theorem}
We apply Theorem~\ref{thmGenSing} to Legendrian singular knot projection with Arnold's contact structure of plane curves to obtain generating (singular) Reidemeister moves.   Note that any first Reidemeister move $\Omega_{1*}$ is forbidden since it changes an equivalence of Legendrian knots.   In this paper, Legendrian fronts from Legendrian knots in the solid torus \cite{ChmutovGoryunov1997} or unit tangent bundle $UTS^2$ \cite{HayanoIto2015} is called a \emph{front presentation} (For \cite{ChmutovGoryunov1997} and \cite{HayanoIto2015}, both possibilities of realizations of Legendrian fronts are essentially the same).  Note also that Legendrian fronts permit only $\Omega_{2c}$ and $\Omega_{2d}$ forbid applying $\Omega_{2a}$ or $\Omega_{2b}$.  
\begin{theorem}\label{thmGenLeg}
We have two statements.  
\begin{enumerate}
\item Generating sets of nonsingular Reidemeister moves are $\Omega_{2c}$ and $\Omega_{2d}$ with any pair: one of $\{ \Omega_{3a},$ $\Omega_{3b}, \Omega_{3c},$ $\Omega_{3d}\}$ and  the other in $\{ \Omega_{3e},$ $\Omega_{3f}, \Omega_{3g},$ $\Omega_{3h} \}$.  
\item Possible front presentations of generating sets of singular Reidemeister moves as the above 
are from the 
triple consisting of $\Omega_{4b}$, $\Omega_{4f}$, and $\Omega_{5d}$.   
\end{enumerate}
Further,  front presentations of $\Omega_{4b}$, $\Omega_{4f}$, and $\Omega_{5d}$ are as in Figures~\ref{ronbunL}, \ref{4f}, and \ref{ronbunN}, respectively.  
\end{theorem}
\begin{figure}
\includegraphics[width=8cm]{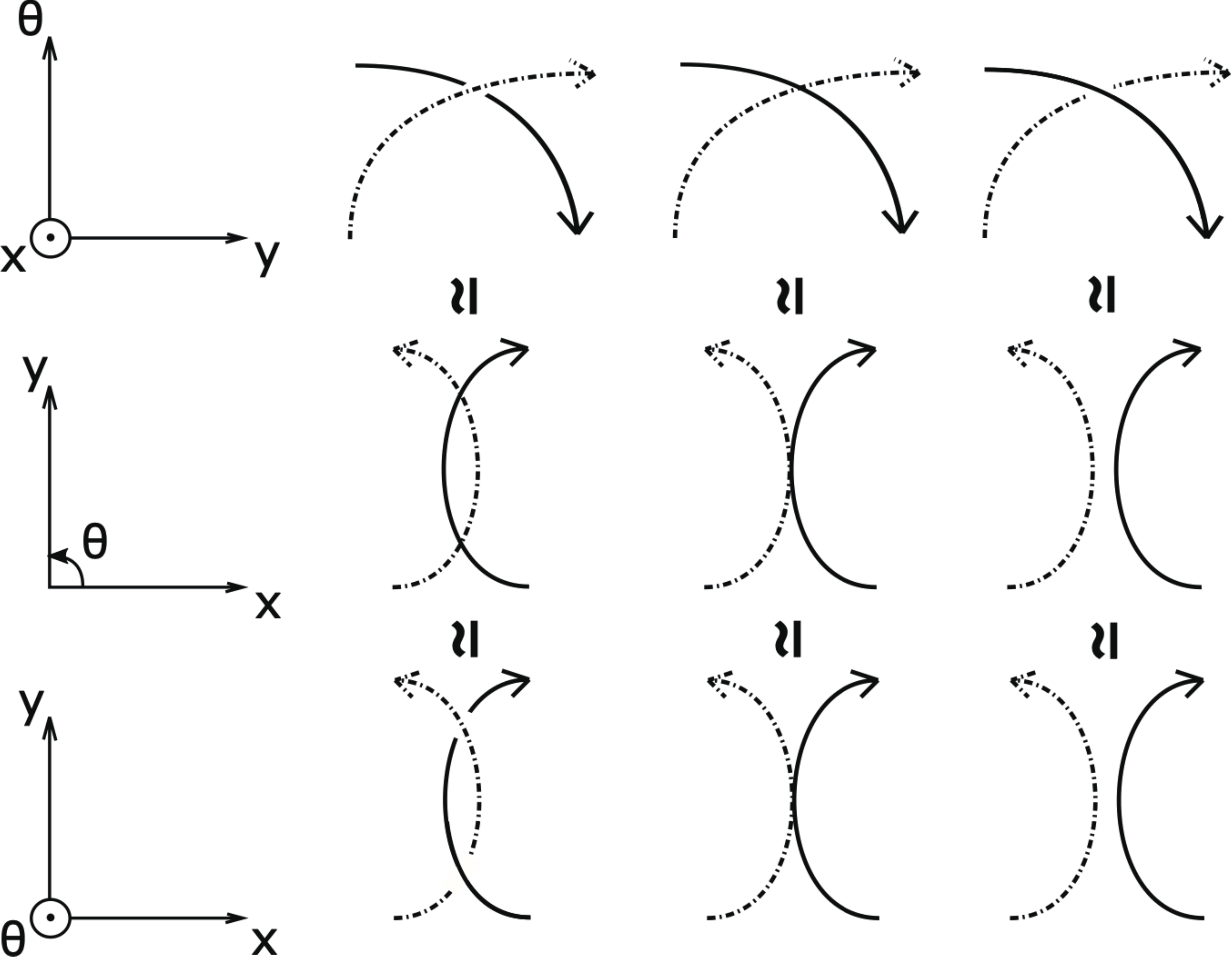}
\caption{Crossing change in the coordinate $(x, y, \theta)$ (1st line) }\label{ronbunK}
\end{figure}

\section{Proof of Theorem~\ref{thmGenSing}}
\cite[Lemmas~5.3--5.8]{BatainehKhaledElhamdadiHajij2018} implies the statement for $\Omega_{4*}$ and \cite[Lemmas~5.10--5.14]{BatainehKhaledElhamdadiHajij2018} implies the statement for $\Omega_{5*}$.  See Tables~\ref{BKEHa} and \ref{BKEHb}.    
\begin{table}[h!]
\caption{}\label{BKEHa}
\begin{tabular}{|c|c|c|c|c|} \hline
$\Omega_{4*}$ & $\Omega_{2*}$ increasing crossings & $\Omega_{4*}$ & $\Omega_{2*}$ decreasing crossings  & Reference  \\ \hline
$\Omega_{4b}$ & $\Omega_{2a}$ & $\Omega_{4d}$  & $\Omega_{2b}$ & \cite[Lemma~5.7]{BatainehKhaledElhamdadiHajij2018} \\ \hline
$\Omega_{4c}$ & $\Omega_{2c}$ & $\Omega_{4a}$  & $\Omega_{2d}$ & \cite[Lemma~5.3]{BatainehKhaledElhamdadiHajij2018} \\ \hline
$\Omega_{4d}$ & $\Omega_{2c}$ & $\Omega_{4a}$  & $\Omega_{2d}$ & \cite[Lemma~5.4]{BatainehKhaledElhamdadiHajij2018} \\ \hline
$\Omega_{4f}$ & $\Omega_{2b}$ & $\Omega_{4h}$  & $\Omega_{2a}$ & \cite[Lemma~5.8]{BatainehKhaledElhamdadiHajij2018} \\ \hline
$\Omega_{4g}$ & $\Omega_{2c}$ & $\Omega_{4e}$  & $\Omega_{2d}$ & \cite[Lemma~5.5]{BatainehKhaledElhamdadiHajij2018} \\   \hline
$\Omega_{4h}$ & $\Omega_{2c}$ & $\Omega_{4e}$  & $\Omega_{2d}$ & \cite[Lemma~5.6]{BatainehKhaledElhamdadiHajij2018} \\ \hline
\end{tabular}
\end{table}

\begin{table}[h!]
\caption{}\label{BKEHb}
\begin{tabular}{|c|c|c|c|c|c|c|c|} \hline
$\Omega_{5*}$ & $\Omega_{1*}$ or $\Omega_{2*} $ & $\Omega_{4*}$ & $\Omega_{5*}$ &  $\Omega_{4*}$ &  $\Omega_{1*}$ or $\Omega_{2*} $   & Reference  \\ \hline
$\Omega_{5b}$ & $\Omega_{1a}$ & $\Omega_{4a}$ & $\Omega_{5d}$  & $\Omega_{4e}$ & $\Omega_{1a}$ & \cite[Lemma~5.10]{BatainehKhaledElhamdadiHajij2018} \\ \hline
$\Omega_{5c}$ & $\Omega_{1b}$ & $\Omega_{4a}$ & $\Omega_{5d}$  & $\Omega_{4e}$ & $\Omega_{1b}$ & \cite[Lemma~5.11]{BatainehKhaledElhamdadiHajij2018} \\ \hline
$\Omega_{5d}$ & $\Omega_{2b} \times 2$ &   & $\Omega_{5a}$ &  & $\Omega_{2a} \times 2$ & \cite[Lemma~5.14]{BatainehKhaledElhamdadiHajij2018} \\ \hline
$\Omega_{5e}$ & $\Omega_{1c}$ & $\Omega_{4e}$  & $\Omega_{5a}$ & $\Omega_{4e}$ & $\Omega_{1a}$ & \cite[Lemma~5.12]{BatainehKhaledElhamdadiHajij2018} \\ \hline
$\Omega_{5f}$ & $\Omega_{1d}$ & $\Omega_{4e}$ & $\Omega_{5a}$  & $\Omega_{4a}$ & $\Omega_{1d}$ & \cite[Lemma~5.13]{BatainehKhaledElhamdadiHajij2018} \\   \hline
\end{tabular}
\end{table}

\section{Proof of Theorem~\ref{thmGenLeg}}
\subsection{Proof for nonsingular Reidemeister moves}
We will show two lemmas: 
\begin{lemma}\label{lem3a}
Any element of $\{\Omega_{3b},  \Omega_{3c}, \Omega_{3d}, \Omega_{3e}, \Omega_{3f}, \Omega_{3g}, \Omega_{3h} \}$ generates $\{ \Omega_{3a}, \Omega_{3h} \}$ using $\Omega_{2c}$ and $\Omega_{2d}$.  
\end{lemma}
See Table~\ref{PolyakTable}.  
\begin{table}[h!]
\caption{Relation between two types of $\Omega_{3}$ via $\Omega_{2c}$ and $\Omega_{2d}$}\label{PolyakTable}
\begin{tabular}{|c|c|c|c|c|} \hline
$\Omega_{3*}$ & $\Omega_{2*}$ increasing crossings & $\Omega_{3*}$ & $\Omega_{2*}$ decreasing crossings  & Reference  \\ \hline
$\Omega_{3b}$ & $\Omega_{2c}$ & $\Omega_{3a}$  & $\Omega_{2d}$ & \cite[Lemma~2.3]{Polyak2010} \\ \hline
$\Omega_{3c}$ & $\Omega_{2c}$ & $\Omega_{3a}$  & $\Omega_{2d}$ & \cite[Lemma~2.4]{Polyak2010} \\ \hline
$\Omega_{3d}$ & $\Omega_{2c}$ & $\Omega_{3a}$  & $\Omega_{2d}$ & Figure~\ref{OhmD} \\ \hline
$\Omega_{3e}$ & $\Omega_{2d}$ & $\Omega_{3h}$  & $\Omega_{2c}$ & Figure~\ref{OhmE} \\ \hline
$\Omega_{3f}$ & $\Omega_{2d}$ & $\Omega_{3a}$  & $\Omega_{2c}$ & \cite[Lemma~2.6]{Polyak2010} \\ \hline
$\Omega_{3g}$ & $\Omega_{2c}$ & $\Omega_{3h}$  & $\Omega_{2d}$ & Figure~\ref{OhmG} \\ \hline
\end{tabular}
\end{table}
\begin{figure}[h!]
\includegraphics[width=10cm]{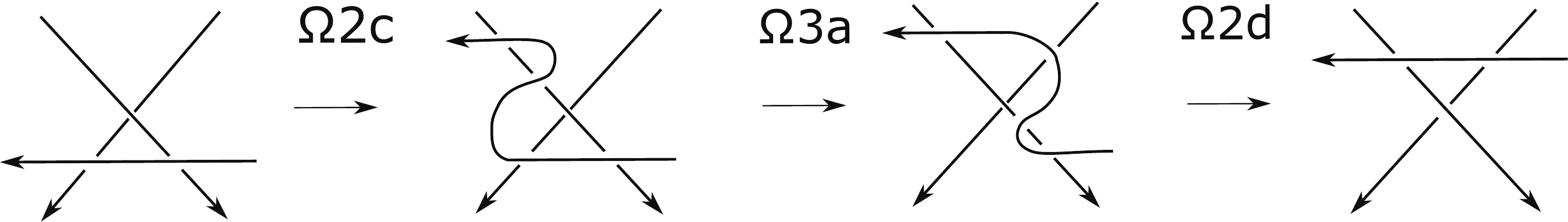}
\caption{$\Omega_{3d}$}\label{OhmD}
\end{figure}
\begin{figure}[h!]
\includegraphics[width=10cm]{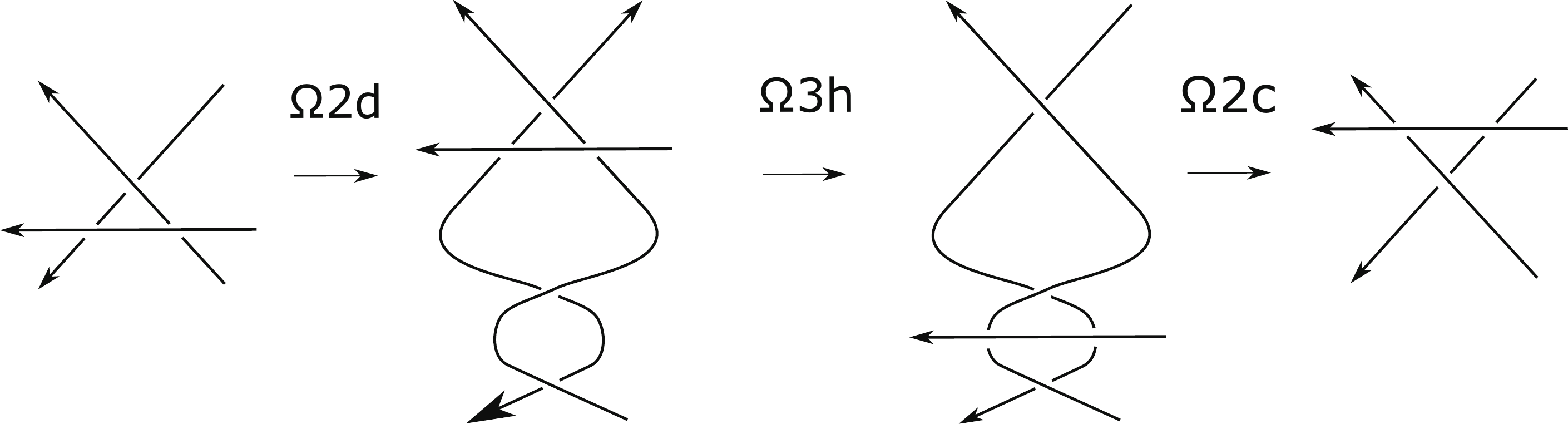}
\end{figure}
\begin{figure}[h!]
\caption{$\Omega_{3e}$}\label{OhmE}
\end{figure}
\begin{figure}[h!]
\includegraphics[width=10cm]{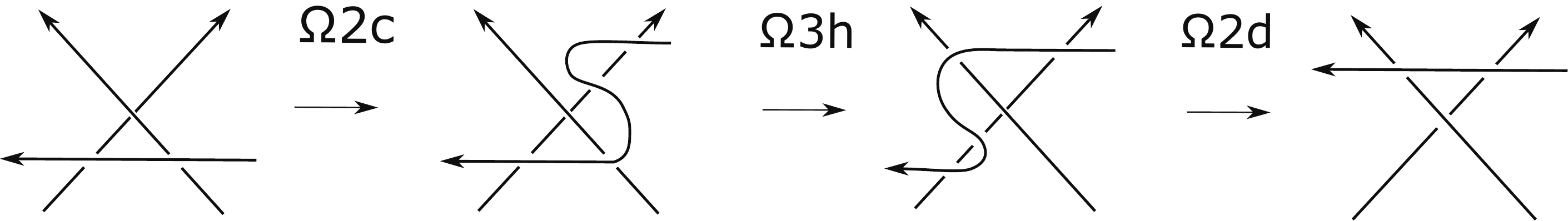}
\caption{$\Omega_{3g}$}\label{OhmG}
\end{figure}
Lemma~\ref{lem3a} directly implies Lemma~\ref{lemPair}.  
\begin{lemma}\label{lemPair}
For any pair of two elements of $\{ \Omega_{3c}, \Omega_{3d}, \Omega_{3e}, \Omega_{3f}, \Omega_{3g}, \Omega_{3h} \}$, one generates another using $\Omega_{2c}$ and $\Omega_{2d}$.   
\end{lemma}
Lemma~\ref{lemPair} implies the proof of the former part of Theorem~\ref{thmGenLeg}.  
\subsection{Proof for singular Reidemeister moves}
For given diagrams in the coordinate the bottom line of Figure~\ref{ronbunK}, which directly deduces the center line, we check the possibility of the realization of the highest line of Figure~\ref{ronbunK} where its realization does not imply any contradiction to the angle $\theta$ of the center line.  
When we try to realize diagrams from $\Omega_{4*}$ or $\Omega_{5*}$ for every case, only the cases as in Figs.~\ref{ronbunL}--\ref{ronbunN} survive.  
The other cases give plane curves locally from $\Omega_{4*}$ or $\Omega_{5*}$, but it implies a contradiction to two tangent lines at self-tangency should be equivalent corresponding to an intersection (that is singular point) of two branches.   

For example, we focus on $\Omega_{4a}$ here.  The main point is the existence of  the same $\theta$ for three branches in the move $\Omega_{4a}$.   Since the singular point is presented by self-tangency (Figure~\ref{ronbunK}), then we have the center column.  However, there is no quadrant the common $\theta$ three branches take. Then we notice that the direction of moving branch should be put between the other two branches, which implies that only $\Omega_{4b}$ and $\Omega_{4f}$ survive.   

For $\Omega_{5*}$, inverse-self-tangency has no common tangency, thus, only $\Omega_{5a}$ and $\Omega_{5d}$ are possible; however, by comparing angles derived from two branches of the (real) crossing, only $\Omega_{5d}$ survives.  
\begin{figure}
\centering
\includegraphics[width=8cm]{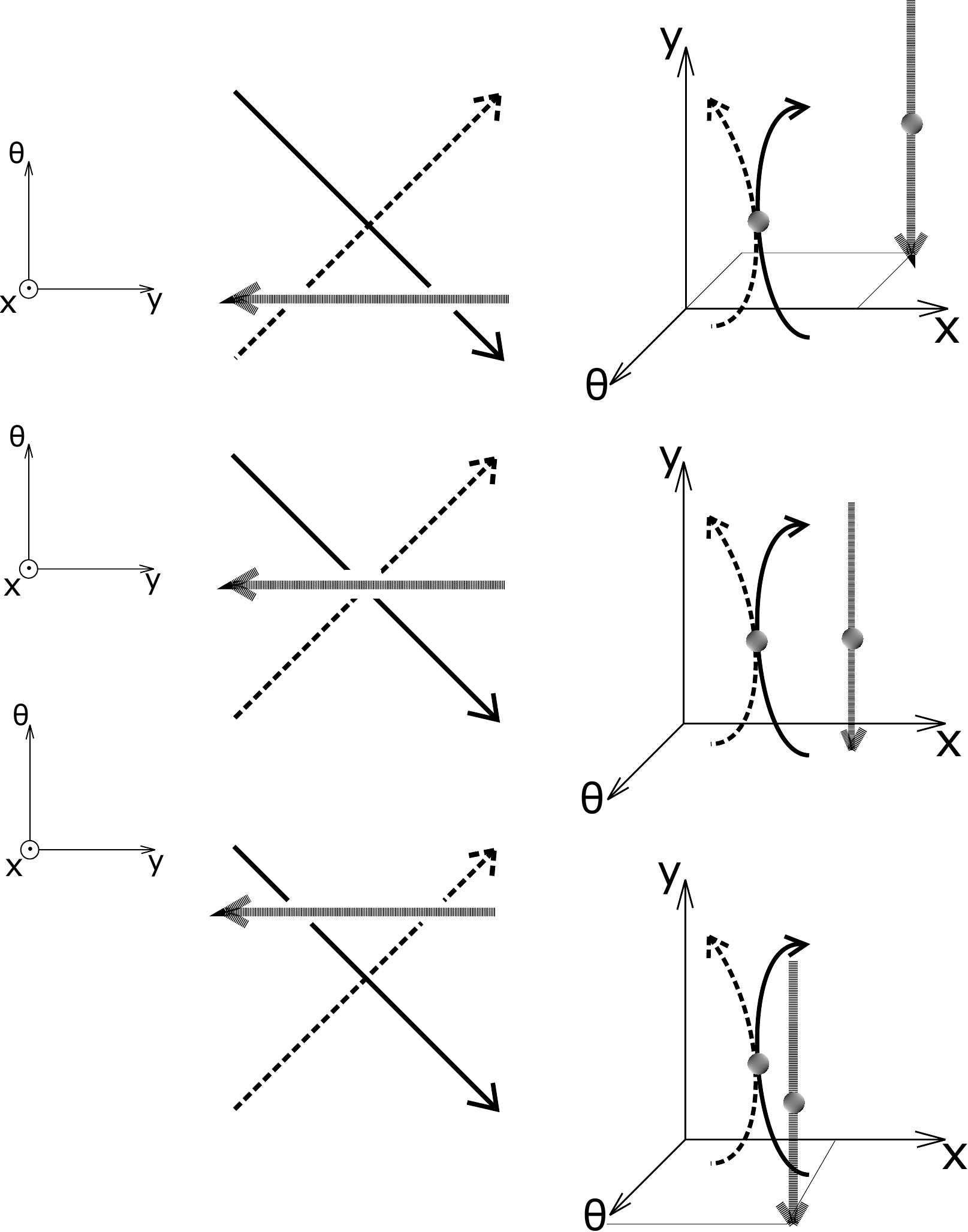}
\caption{Front presentation of $\Omega_{4a}$.  }\label{ronbunM}
\end{figure}

\section*{Acknowledgements}
The work of the authors is partially supported by JSPS KAENHI Grant  number 20K03604 
and by the alumni of National Institute of Technology, Ibaraki College.     

\bibliographystyle{plain}
\bibliography{Ref}
\end{document}